\documentclass[runningheads]{llncs}
\usepackage{graphicx}

\usepackage{amsfonts}
\usepackage{amsmath}

\usepackage{caption}
\usepackage{subfigure}

\usepackage{booktabs}
\usepackage{float}
\usepackage{color}
\usepackage{multirow}

\newcommand{\floor}[1]{\left\lfloor #1 \right\rfloor}

\DeclareUnicodeCharacter{2212}{-}

\begin{document}
\title{Tightening Discretization-based MILP Models for the Pooling Problem using Upper Bounds on Bilinear Terms}
\titlerunning{Tightening Discretization-based MILP Models for Pooling}

%
\author{Yifu Chen\inst{1} \and
Christos T. Maravelias\inst{2,3,*} \and Xiaomin Zhang\inst{4}}
\authorrunning{Chen et al.}
%
\institute{Department of Chemical and Biological Engineering, University of Wisconsin-Madison, Madison, WI 53706, USA \and
Department of Chemical and Biological Engineering, Princeton University, Princeton, NJ 08544, USA \and
Andlinger Center for Energy and the Environment, Princeton University, Princeton, NJ 08544, USA \and
Department of Computer Sciences, University of Wisconsin-Madison, Madison, WI 53706, USA\\
\email{*Corresponding author: maravelias@princeton.edu}}
\maketitle              
\begin{abstract}
Discretization-based methods have been proposed for solving nonconvex optimization problems with bilinear terms such as the pooling problem. These methods convert the original nonconvex optimization problems into mixed-integer linear programs (MILPs). In this paper we study tightening methods for these MILP models for the pooling problem, and derive valid constraints using upper bounds on bilinear terms. Computational results demonstrate the effectiveness of our methods in terms of
reducing solution time.
\end{abstract}

\section{Introduction}\label{Introduction}

The pooling problem is a nonconvex optimization problem that has been studied extensively. First studied by Harvey \cite{Haverly1978}, it continues to be an active research topic. In the pooling problem, we are given a set of streams to be blended in a set of pools to produce certain products. Streams, pools, and products are connected through pipelines. Streams possess multiple properties and products have specifications on those properties. We aim to optimize the flows through the pipelines to maximize profit, subject to product specifications and other constraints, such as pool and pipeline capacity.

Various formulations for the pooling problem have been proposed (\cite{Haverly1978,Ben-Tal1994,Tawarmalani2002,Alfaki2013,Boland2016}), and a number of extensions of the pooling problem have been studied. Comprehensive reviews on these topics can be found in a survey by Misener and Floudas \cite{Misener2009} and a thesis by Gupte \cite{Gupte2012}. 

Solving the pooling problem is computationally challenging \cite{Haugland2016}, and many solution methods for the pooling problem have been proposed. These solution methods include, for example, methods based on piecewise linear relaxations (\cite{Gounaris2009,Castro2015}) and domain partitioning (\cite{Misener2011,Nagarajan2019}), as well as methods that exploit the product specifications constraints (\cite{Papageorgiou2012,Chen2020}). Nonlinear relaxations of the pooling problem have also been studied. For example, Kimizuka et al. studied the second order cone relaxation of the pooling problem \cite{Kimizuka2019}, Luedtke et al. proposed a strong convex nonlinear relaxation derived from extended formulation \cite{Luedtke2020}, and Dey et al. studied the convexifiation methods for rank-one-based substructures arising in the pooling problem \cite{Dey2020}.  

Discretization-based methods for solving the pooling problem have received much attention over the years. Since the only nonlinear component in the pooling problem are bilinear terms that contain two continuous variables, by discretizing one variable in each bilinear term, one can convert the original nonconvex optimization problem into a mixed-integer linear program (MILP). Kolodziej et al. proposed global optimization algorithms based on solving MILPs resulted from discretization \cite{Kolodziej2013}. Dey and Gupte analyzed MILP techniques to address bilinear terms \cite{Dey2015}, Gupte et al. studied relaxations and discretizations for the pooling problem \cite{Gupte2017}. Ceccon and Misener \cite{Ceccon} implemented a discretization technique into an extensible solver for the pooling problem.  

In this paper, we study tightening constraints for the discretization-based MILP models for the pooling problem, with the focus on using the upper bounds on bilinear terms to tighten the MILP models. Specifically, we consider the following set, which, as we will show later, arises in the models for the pooling problem:
\begin{equation}\label{Def S_1}
\textbf{S}_1 = \{(x,y) \in \mathbb{R}_{+}^{2}:  x \leq x^{\text{U}}, y \leq 1, xy \leq w^{\text{U}}\}
\end{equation}
Note that setting the upper bound on $y$ to one is not overly restrictive since one can re-scale variables $x$ and $y$ accordingly. The upper bound on the bilinear term, $w^{\text{U}}$, is said to be nontrivial if $w^{\text{U}} < x^{\text{U}}$. In this paper, we assume such an upper bound is nontrivial.

We start with discretizing $y$. Specifically, we restrict the value of $y$ to be in a discrete set $\textbf{O} = \{0,2^{1-n}, 2^{2 - n}, ... ,1 - 2^{1 - n}, 1\}$, whose elements form an arithmetic sequence, where $n$ is a positive integer. After discretizing $y$, from set $\textbf{S}_1$ we have:

\begin{equation}\label{Def_S_1_star}
\textbf{S}_1^{*} = {\{(x,y) \in \mathbb{R}_{+}^{2}: x \leq x^{\text{U}}, y \in \textbf{O}, xy \leq w^{\text{U}}\}}
\end{equation}

We note that $\textbf{S}_1^{*}$ can be represented with linear constraints. First, to model $y \in \textbf{O}$, we consider the binary expansion of $y$:

\begin{equation}\label{Binary_expansion}
y= \sum_{p \in [\![n]\!]} 2^{1-p}Z_{p}
\end{equation}
where $Z_p \in \{0,1\}$ and $[\![n]\!] = \{1,2,...,n\}$. By introducing binary variable $Z_p$, \eqref{Binary_expansion}, and a continuous variable $w$, we have the following set:
\begin{equation}\label{Def_S_1_prime}
\textbf{S}_1^{'} = {\{(w,x,y,Z) \in \mathbb{R}_{+}^{3} \times \{0,1\}^{n}: w \leq w^{\text{U}}, x \leq x^{\text{U}}, y \leq 1, \eqref{Binary_expansion}, w=xy\}}
\end{equation}
We next re-write $w=xy$ as:

\begin{equation}\label{Re-Def_w}
w=\sum_{p } 2^{1-p}xZ_{p}
\end{equation}

By defining $v_p = xZ_{p}$ for each $p \in [\![n]\!]$, we can linearize \eqref{Re-Def_w} exactly using:

\begin{equation}\label{Re-Def_w_v}
w=\sum_{p } 2^{1-p}v_{p}
\end{equation}

\begin{equation}\label{UB_v_p_1}
v_{p} \leq x^{\text{U}} Z_{p}, \quad \forall p 
\end{equation}

\begin{equation}\label{UB_v_p_2}
v_{p} \leq x, \quad \forall p 
\end{equation}

\begin{equation}\label{LB_v_p_1}
v_{p} \geq x + x^{\text{U}} Z_{p} - x^{\text{U}}, \quad \forall p 
\end{equation}

We define the set $\textbf{S}_2$ as: 
\begin{multline}\label{Def S_2}
\textbf{S}_2 = \{(v,w,x,y,Z) \in \mathbb{R}_{+}^{n} \times \mathbb{R}_{+}^{3} \times \{0,1\}^{n}: \\
w \leq w^{\text{U}}, x \leq x^{\text{U}}, y \leq 1, \eqref{Binary_expansion}, \eqref{Re-Def_w_v} - \eqref{LB_v_p_1}\}
\end{multline}
which is a linear representation of $\textbf{S}_1^{*}$. 


We note one related work from Gupte et al.\cite{Gupte2013}, in which the authors studied a set equivalent to
\begin{equation*}
\textbf{S}_{2}^{'} = \{(v,w,x,y,Z) \in \mathbb{R}_{+}^{n} \times \mathbb{R}_{+}^{3} \times \{0,1\}^{n}: x \leq x^{\text{U}}, y \leq 1, \eqref{Binary_expansion},\eqref{Re-Def_w_v} - \eqref{LB_v_p_1}\}
\end{equation*}

which is a relaxation of $\textbf{S}_2$ since it can be obtained by relaxing the nontrivial upper bound on $w$ (see (2.4) in \cite{Gupte2013}). The authors presented the convex hull of $\textbf{S}_{2}^{'}$ by exploiting the knapsack structure resulted from the binary expansion. 

Different from \cite{Gupte2013}, in this paper we study valid constraints for $\textbf{S}_2$ derived from nontrivial upper bounds on bilinear terms. Specifically, we give the full description of the convex hull of $\textbf{S}_1^{*}$; in addition, we tighten \eqref{UB_v_p_1} using tightened upper bounds on $x$. We also note that while our methods are presented in the context of the pooling problem, they can be extended to address other nonconvex optimization problems with nontrivial upper bounds on bilinear terms.

We note that methods that exploit the nontrivial bounds of bilinear terms have been proposed; most notably, Anstreicher et al. derived convex hull representations for bounded products of continuous variables
 \cite{Anstreicher2021}. In this work we focus on the discretized version of the problem, and propose a different set of tightening constraints.  

The remainder of the paper is organized as follows. In Section \ref{Background} we present background material, including the problem statement, a nonlinear model for the pooling problem, and an MILP model obtained from discretization. In Section \ref{Valid Constraints} we derive valid constraints for the MILP model that utilizes bounds on bilinear terms, and present numerical examples. In Section \ref{Computational Results}, we demonstrate the effectiveness of our constraints, showing computational results for models derived from different formulations for (1) the pooling problem and (2) discretization-based pooling models, as well as comparisons with other methods.

\section{Background}\label{Background}
\subsection{Problem statement}
In the standard setting, the pooling problem is defined in terms of the indices/sets and parameters given in Table \ref{tab:Sets and parameters}:

\begin{table}[H]
\centering
\captionsetup{font = normalsize}
\caption {Sets and parameters for the pooling problem} \label{tab:Sets and parameters} 
\begin{tabular}{ll}
\toprule
\textit{Indices/Sets}: & \\
$   i \in $ \textbf{I}:   & Inputs (streams) \\
$   j \in $ \textbf{J}: & Pools \\
$   k \in $ \textbf{K}: & Products \\
$   l \in $ \textbf{L}: & Properties \\
\textit{Parameters}: & \\
$\alpha_i $:   & Unit cost of stream $i$ \\
$\beta_k $:   & Price of product  $k$ \\
$\gamma_j $: & Capacity of pool $j$ \\
$\upsilon_{jk}  $: & Capacity of the pipeline between pool $j$ and product $k$ \\
$\pi_{il} $: & Value of property $l$ for stream $i$ \\
$\psi_{kl} $: & Upper bounding specification for property $l$ for product $k$ \\
$\omega_{k} $: & Maximum demand for product $k$ \\ \bottomrule
\end{tabular}
\end{table}

For any product, the combined flows from all pools to that product must satisfy the corresponding specifications. We aim to find flows (from streams to pools and from pools to products) that maximize profit. We assume that there are no flows between pools, no stream flow accumulation in pools, and all product properties are the average of the properties of the streams blended weighted by volume fraction. Without loss of generality, we assume we have only upper bounding specifications on the properties.

\subsection{A nonlinear model for the pooling problem}
\label{sec:A nonlinear model}
Various models have been proposed for the pooling problem. In this paper, we first study a model that contains variable for split fraction, which is similar to the one presented in Boland et al. \cite{Boland2016} and the source-based model studied by Lotero et al. \cite{Lotero2016}. We define the following nonnegative continuous variables: \newline

\begin{tabular}{ll}
$F_{ij} $:   & Flow of stream $i$ to pool $j$ \\
$\hat{F}_{ijk} $:   & Flow of stream $i$ from pool $j$ to product $k$  \\
$R_{jk} $:   & Split fraction for total inlet flows from pool $j$ to product $k$ ($R_{jk} \leq 1$)  \\
\end{tabular} \newline

We have the following constraints: \newline

Pool capacity:\newline
\begin{equation}\label{Pool capacity}
\sum_{i  } F_{ij} \leq \gamma_j,  \quad   \forall j  
\end{equation}

Product demand:\newline
\begin{equation}\label{Product demand}
\sum_{i  }\sum_{j  } \hat{F}_{ijk} \leq \omega_{k},  \quad   \forall k  
\end{equation}

Product specifications:\newline
\begin{equation}\label{Product specifications}
\sum_{i  }\sum_{j  } \pi_{il}\hat{F}_{ijk} \leq \psi_{kl}\sum_{i  }\sum_{j  } \hat{F}_{ijk},  \quad  \forall k  , l  
\end{equation}

Stream splitting:\newline
\begin{equation}\label{Stream splitting}
\hat{F}_{ijk} = F_{ij}R_{jk},  \quad   \forall i  , j  , k  
\end{equation}

Note that \eqref{Stream splitting} is an equality constraint with a bilinear term.

For split fraction $R_{jk}$ we have:\newline
\begin{equation}\label{RSimplex}
\sum_{k  } R_{jk} = 1,  \quad   \forall j  
\end{equation}

\eqref{Stream splitting} and \eqref{RSimplex} enforce that there is no flow accumulation in pools.

We also have upper bound on the flows from pools to products due to pipeline capacity:\newline
\begin{equation}\label{Upper bound on the flows from pools to products}
\sum_{i  } \hat{F}_{ijk} \leq \upsilon_{jk},  \quad   \forall j  , k  
\end{equation}

We note that the left hand side (LHS) of \eqref{Upper bound on the flows from pools to products} is also upper bounded by $\gamma_j$. However, $\upsilon_{jk}$ is a tighter upper bound since in practice pipeline capacity is typically much smaller than pool capacity. In this work we assume $\upsilon_{jk} < \gamma_j$. In other words, $\upsilon_{jk}$ is a nontrivial upper bound on the bilinear term $\sum_{i} \hat{F}_{ijk} = \sum_{i  } F_{ij}R_{jk}$.

Reformulation–Linearization Technique (RLT) constraints can be added to tighten the formulation. Summing over index $k$ on both sides of \eqref{Stream splitting}, we have:

\begin{equation*}
\sum_{k  } \hat{F}_{ijk} = F_{ij} \sum_{k  } R_{jk},  \quad   \forall i  , j  
\end{equation*}

which, combined with \eqref{RSimplex}, leads to:

\begin{equation}\label{Sum_k F_hat}
\sum_{k  } \hat{F}_{ijk} = F_{ij},  \quad  \forall i  , j  
\end{equation}

Another RLT constraint can be obtained by multiplying both sides of \eqref{Pool capacity} with $R_{jk}$ (a nonnegative variable):

\begin{equation*}
\sum_{i  } F_{ij}R_{jk} \leq \gamma_jR_{jk},  \quad \forall  j  , k  
\end{equation*}

which, combined with \eqref{Stream splitting}, leads to:

\begin{equation}\label{Sum_i F_hat}
\sum_{i  } \hat{F}_{ijk} \leq \gamma_jR_{jk},  \quad  \forall j  , k  
\end{equation}

The objective function is profit maximization:

\begin{equation}\label{Obj: Max profit}
\text{max} \sum_{i  } \sum_{j  } (\sum_{k  } \beta_k \hat{F}_{ijk} - \alpha_i {F}_{ij}  )
\end{equation}

\eqref{Pool capacity} - \eqref{Obj: Max profit} comprise a nonlinear model for the pooling problem. We next present a discretization-based model derived from it, and show how we exploit the structure of the above nonlinear model to tighten the discretization-based model. 

\subsection{A discretization-based model for the pooling problem}
\label{sec:A discretization-based model}

We consider a discretization-based model for the above nonlinear model by restricting $R_{jk} \in \textbf{O}$, with the following binary expansion of $R_{jk}$:

\begin{equation}\label{Binary expansion R}
R_{jk} = \sum_{p \in [\![n]\!]  } 2^{1-p}Z_{jkp}, \quad  \forall j  , k  
\end{equation}

Combining \eqref{Stream splitting} and \eqref{Binary expansion R} we have:

\begin{equation}\label{expansion F hat quadratic}
\hat{F}_{ijk} = \sum_{p  } 2^{1-p}{F}_{ij}Z_{jkp}, \quad  \forall i  , j  , k  
\end{equation}

We introduce a nonnegative continuous variable {$V_{ijkp}$} to represent {${F}_{ij}Z_{jkp}$}:

\begin{equation}\label{Def V}
V_{ijkp} = {F}_{ij}Z_{jkp}, \quad  \forall i  , j  , k  , p  
\end{equation}

Since ${F}_{ij}$ is upper bounded by $\gamma_j$, we linearize \eqref{expansion F hat quadratic} using constraints similar to \eqref{Re-Def_w_v} - \eqref{LB_v_p_1}:

\begin{equation}\label{Expansion F hat linear}
\hat{F}_{ijk} = \sum_{p  } 2^{1-p}V_{ijkp}, \quad  \forall i  , j  , k  
\end{equation}

\begin{equation}\label{Linearizing V 1}
V_{ijkp} \leq \gamma_jZ_{jkp}, \quad \forall  i  , j  , k  , p  
\end{equation}

\begin{equation}\label{Linearizing V 2}
V_{ijkp} \leq {F}_{ij}, \quad \forall  i  , j  , k  , p  
\end{equation}

\begin{equation}\label{Linearizing V 3}
V_{ijkp} \geq {F}_{ij} - \gamma_j(1 - Z_{jkp}), \quad \forall  i  , j  , k  , p  
\end{equation}

Inequalities (\ref{Linearizing V 1}) and (\ref{Linearizing V 3}) can be tightened. To simplify notation, we define $\tilde{F}_{j} = \sum_{i} F_{ij}$, and $\Bar{V}_{jkp} = \sum_{i} {V}_{ijkp}$. By summing over index $i$ for \eqref{Def V} we have: 

\begin{equation}\label{Def V Agg}
\Bar{V}_{jkp} = \tilde{F}_{j}Z_{jkp}, \quad  \forall j  , k  , p  
\end{equation}

From \eqref{Pool capacity} recall that $\gamma_j$ is also the upper bound on $\tilde{F}_{j}$. Hence, we have the following two tightening constraints:

\begin{equation}\label{Linearizing V 1 Agg}
\Bar{V}_{jkp} \leq \gamma_jZ_{jkp}, \quad \forall  j  , k  , p  
\end{equation}

\begin{equation}\label{Linearizing V 3 Agg}
\Bar{V}_{jkp} \geq \tilde{F}_{j} - \gamma_j(1 - Z_{jkp}), \quad  \forall j  , k  , p  
\end{equation}

\eqref{Pool capacity} - \eqref{Product specifications}, \eqref{RSimplex} - \eqref{Binary expansion R}, \eqref{Expansion F hat linear} - \eqref{Linearizing V 3}, and \eqref{Linearizing V 1 Agg} - \eqref{Linearizing V 3 Agg} comprise a discretization-based model for the pooling problem, henceforth referred to as $\text{M}^{\text{SB}}$. Solving $\text{M}^{\text{SB}}$ leads to a feasible solution to the original nonlinear model of the pooling problem, and the optimal objective function value of $\text{M}^{\text{SB}}$ is a lower bound on the optimal objective function value of the nonlinear problem. In practice, discretization-based models provide pretty good solutions within a reasonable amount of time; this observation is also noted by Dey and Gupte after extensive computational comparisons \cite{Dey2015}.

$\text{M}^{\text{SB}}$ contains a structure similar to the set $\textbf{S}_2$ presented in Section \ref{Introduction}; specifically, for a given $(j,k)$ pair, one can consider $\sum_{i} \hat{F}_{ijk}$ as $w$, $\sum_{i} {F}_{ij}$ as $x$, $R_{jk}$ as $y$, and $\sum_{i} {V}_{ijkp}$ as $v_p$, while \eqref{Upper bound on the flows from pools to products}, \eqref{Expansion F hat linear}, \eqref{Linearizing V 2}, \eqref{Linearizing V 1 Agg} - \eqref{Linearizing V 3 Agg} are similar to the constraints involved in $\textbf{S}_2$.  We next present tightening constraints for $\text{M}^{\text{SB}}$, with particular interest in constraints that are associated with flows into pools.

We note that many other formulations for the pooling problem have been proposed (\cite{Haverly1978,Ben-Tal1994,Tawarmalani2002,Alfaki2013,Boland2016}), and we show results for our methods applied to another formulation in Section \ref{Computational Results}.

\section{Valid Constraints}\label{Valid Constraints}

We first provide a high level description of the structure in $\text{M}^{\text{SB}}$ that we exploit to generate valid constraints. Let $\Bar{F}_{jk} = \sum_{i} \hat{F}_{ijk}$ and consider the bilinear term $\Bar{F}_{jk} = \tilde{F}_{j}R_{jk}$. For a given $(j,k)$ pair, we define the following set:
\begin{equation}\label{Def_S_jk}
\textbf{S}_{jk} = \{(R_{jk},\tilde{F}_{j},\Bar{F}_{jk}) \in \mathbb{R}_{+}^{3}:  \tilde{F}_{j} \leq \gamma_j, \Bar{F}_{jk} \leq \upsilon_{jk}, R_{jk} \leq 1, \Bar{F}_{jk} = \tilde{F}_{j}R_{jk}   \}
\end{equation}

and after restricting $R_{jk} \in \textbf{O}$, we have the following set:

\begin{equation}\label{Def_S_jk_star}
\textbf{S}_{jk}^{*} = \{(R_{jk},\tilde{F}_{j}) \in \mathbb{R}_{+}^{2}:  \tilde{F}_{j} \leq \gamma_j, \tilde{F}_{j}R_{jk} \leq \upsilon_{jk}, R_{jk} \in \textbf{O}   \}
\end{equation}

In Section \ref{Rounding cuts}, we present a set of valid constraints derived from $\textbf{S}_{jk}^{*}$. Specifically, we analyze the feasible space defined by it and provide a closed-form representation of its convex hull.

Similar to $\textbf{S}_{1}^{*}$, $\textbf{S}_{jk}^{*}$ can be represented by linear constraints. Such a linear representation requires linearizing $\tilde{F}_jZ_{jkp}$, which utilizes the upper bound on $\tilde{F}_{j}$. We note that while $\gamma_j$ is a valid upper bound, a tighter upper bound can be obtained in certain situations. Details are presented in Section \ref{p-dependent}. 

For the remaining parts of this section, we focus on a $(j,k)$ pair, and for simplicity we drop these two indices.

\subsection{Rounding cuts for flows into pools}\label{Rounding cuts}

We first analyze the feasible space defined by $\textbf{S}^{*}$:

\begin{equation*}
\textbf{S}^{*} = \{(R,\tilde{F}) \in \mathbb{R}_{+}^{2}:  \tilde{F} \leq \gamma, \tilde{F}R \leq \upsilon, R \in \textbf{O}   \}
\end{equation*}

For a fixed $R$, we have $\tilde{F} \in [0, \text{min}(\gamma, \upsilon/R)]$, a vertical line segment in ($R,\tilde{F}$) space. Since $R$ can only take discrete values in the finite set $\textbf{O}$, the feasible space defined by $\textbf{S}^{*}$ is the union of those vertical line segments. One can thus conclude that the convex hull of $\textbf{S}^{*}$, conv($\textbf{S}^{*}$), is a two-dimensional polyhedron. 

\subsubsection{Illustrative examples} 

\begin{figure}[!htb]
\includegraphics[width=\textwidth]{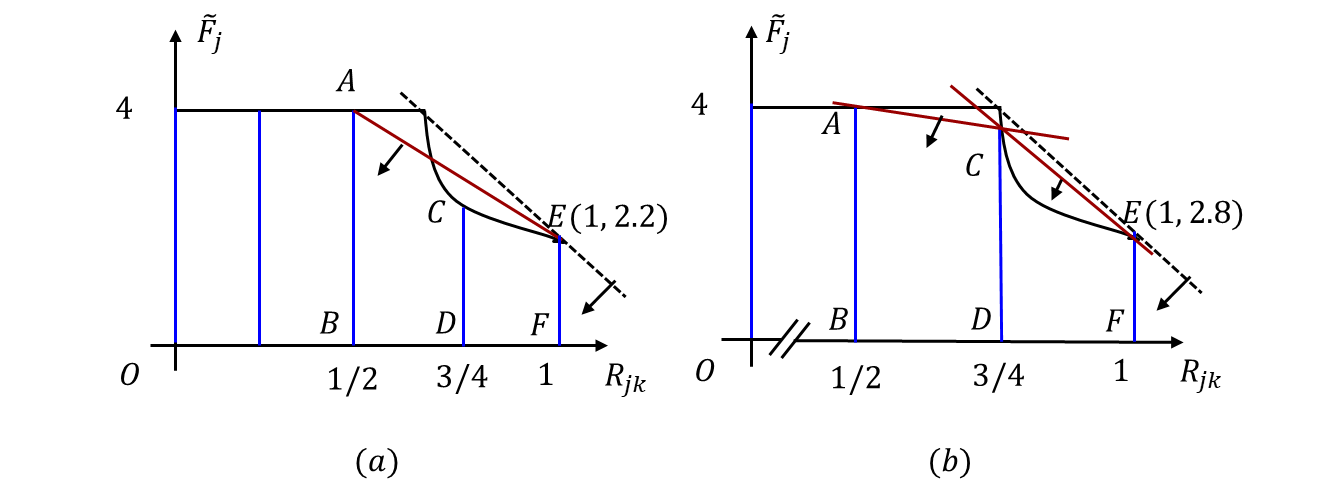}
\caption{Illustrative graph for the examples (a) $\upsilon = 2.2$.  (b) $\upsilon = 2.8$.} \label{IllustrativeFig}
\end{figure}

Consider an example with $\gamma =4$, $\upsilon =2.2$, and $p \in \{1,2,3\}$. We show the feasible space for variables $\tilde{F}$ and $R$ defined by $\textbf{S}^{*}$ in Fig.\ref{IllustrativeFig}(a) (solid blue lines). Here, line segment $AB$ represents the largest value of $R$ for which $\tilde{F}$ can still take value $\gamma$, and line segment $EF$ represents the possible values of $\tilde{F}$ when $R = 1$. It is clear that the constraint represented by the red solid line containing points $A$ and $E$ is a valid constraint for $\textbf{S}^{*}$. We note that at first glance, such a constraint might look similar to the well-known mixed-integer rounding cut \cite{Marchand2001}; we show the difference between the two in Appendix A. 

Consider another example, where everything remains unchanged, except that we now have $\upsilon =2.8$, shown in Fig.\ref{IllustrativeFig}(b). If we again derive a constraint based on points $A$ and $E$, the resulting constraint will not be valid since it will cut off part of the line segment $CD$, which represents the smallest value of $R$ for which $\gamma > \upsilon/R$. Here, we propose two valid constraints represented by red solid lines. We note that for both examples, the proposed constraints are tighter than the constraint included in the linearization (dotted line), which does not utilize the information from the discretization. Furthermore, we note that for both examples, the mentioned constraints, along with bounds on $\tilde{F}$ and $R$, fully describe conv($\textbf{S}^{*}$).

\subsubsection{Valid constraints and polyhedral results}

We next formalize the valid constraints for $\textbf{S}^{*}$. Let $\mu^{-}$ denote the maximum value of $R$ for which $\tilde{F}$ can still take value $\gamma$, we have $\mu^{-}=2^{1-n}\floor{\upsilon/(2^{1-n}\gamma)}$. Further, let $\mu^{+}$ denote the adjacent possible value for $R$ that is greater than $\mu^{-}$, we have $\mu^{+}=2^{1-n}(\floor{\upsilon/(2^{1-n}\gamma)}+1) = \mu^{-} + 2^{1-n}$. We further define the following two parameters: $\delta = (\gamma - \upsilon)/(\mu^{-} - 1)$ and $\epsilon = (\gamma - \upsilon/\mu^{+})/(\mu^{-} - \mu^{+})$ We have the following constraints:

\begin{equation}\label{Rounded Cut 1}
\tilde{F} \leq \delta(R - 1) + \upsilon, \quad \delta \geq \epsilon
\end{equation}

\begin{equation}\label{Rounded Cut 2}
\tilde{F} \leq \epsilon(R - \mu^{-}) + \gamma,\quad \delta < \epsilon
\end{equation}

\begin{equation}\label{Rounded Cut 3}
\tilde{F} \leq \frac{\upsilon/\mu^{+} -\upsilon}{\mu^{+} - 1}(R - 1) + \upsilon, \quad \delta < \epsilon
\end{equation}

We revisit the illustrative examples above. For the example with $\upsilon =2.2$ shown in Fig.\ref{IllustrativeFig}(a), one can verify that $\delta = (4-2.2)/(-0.5)=-3.6$ and $ \epsilon = (4-2.2/0.75)/(-0.25) \newline \approx -4.27$, which represent the slope of line segment $AE$ and $AC$, respectively. Since $\delta > \epsilon$, \eqref{Rounded Cut 1} will be enforced for this example, represented by the red line containing $AE$. For the example with $\upsilon =2.8$ shown in Fig.\ref{IllustrativeFig}(b), we have $\delta < \epsilon$ (since the slope of line segment $AE$ is smaller than $AC$), thus we enforce \eqref{Rounded Cut 2} - \eqref{Rounded Cut 3} for this example, represented by two red lines containing $AC$ and $CE$, respectively. We mention the special case $\mu^{+} = 1$ since ($\mu^{+} - 1$) shows up in the denominator in \eqref{Rounded Cut 3}; we note that division by zero will not occur since when $\mu^{+} = 1$, we have $\delta = \epsilon$ and thus \eqref{Rounded Cut 3} will not be enforced.

Proposition 1 below states that the above three constrains, along with $\tilde{F} \leq \gamma$, $R \leq 1$, describe conv($\textbf{S}^{*}$).

\begin{proposition}\label{P1}conv($\textbf{S}^{*}$) $=\{(R,\tilde{F}) \in \mathbb{R}_{+}^{2}:  \tilde{F} \leq \gamma, R \leq 1, \eqref{Rounded Cut 1} - \eqref{Rounded Cut 3} \}$.\end{proposition}

The full proof of Proposition 1 is presented in Appendix B. Here, we present a sketch of the proof, focusing on the case where $\delta \geq \epsilon$.

We define the polyhedron
\begin{equation*}
Q_1 = \{(R,\tilde{F}) \in \mathbb{R}_{+}^{2}:  \tilde{F} \leq \gamma, R \leq 1, \eqref{Rounded Cut 1} \}
\end{equation*}

Recall that by definition conv($\textbf{S}^{*}$) is the set of all convex combinations of points in $\textbf{S}^{*}$. We next show that (1) $Q_1$ contains all convex combinations of points in $\textbf{S}^{*}$, and (2) all points in $Q_1$ are convex combinations of points in $\textbf{S}^{*}$.

For (1), we first note that since all points in $\textbf{S}^{*}$ are in (finitely many) line segments, all convex combinations of those points can be represented by convex combinations of the end points for those line segments. We show that all these end points are in $Q_1$. Since $Q_1$ is a polyhedron, all convex combinations of these points are also in $Q_1$, thus $Q_1$ contains all convex combinations of points in $\textbf{S}^{*}$.

For (2), one can verify that $Q_1$ contains five extreme points: $P_1 = (0,0),P_2 = (0,\gamma),P_3 = (\mu^{-},\gamma), P_4 = (1,\upsilon), P_5 = (1,0).$ Since all points in $Q_1$ are convex combinations of these extreme points, which are also in $\textbf{S}^{*}$, we conclude that all points in $Q_1$ are convex combinations of points in $\textbf{S}^{*}$.

We show the full proof in Appendix B, in which more details are provided and the case where $\delta < \epsilon$ is also discussed.

\subsection{p-dependent upper bounds on flows into pools}\label{p-dependent}

In this section we focus on the linear representation of the set $\textbf{S}^{*}$, which involves linearizing \eqref{Def V Agg}. Specifically, we note that when certain $Z_{p}$ takes the value 1, the upper bound on $\tilde{F}$ can be tightened, which leads to constraints tighter than \eqref{Linearizing V 1 Agg}.

We first consider the case where for a given $(j,k)$ pair,  {${Z}_{1} = 1$}. In this case, from \eqref{Binary expansion R} we have {${R} = 1$}; in other words, all inlet flows to pool $j$ go to product $k$. Thus, we have $\tilde{F} \leq \upsilon$. We can tighten \eqref{Linearizing V 1 Agg} using:

\begin{equation}\label{Linearizing V 1 tightened 1}
\Bar{V}_{1} \leq \upsilon Z_{1}  
\end{equation}

Now, consider the case where for a given $(j,k)$ pair and $p = 2$,  {${Z}_{2} = 1$}. In this case, from \eqref{Binary expansion R} we know that {${R} \geq 1/2$}, in other words, no less than half the inlet flows to pool $j$ go to product $k$. Thus, we have {$\tilde{F} \leq 2\upsilon$}. We can tighten \eqref{Linearizing V 1 Agg} using:

\begin{equation}\label{Linearizing V 1 tightened 2}
\Bar{V}_{2} \leq 2\upsilon Z_{2}
\end{equation}

Following the same approach we can derive constraints for $p=3,4...$, and the constraints obtained from this approach will be tighter than \eqref{Linearizing V 1 Agg} when $2^{p-1}\upsilon < \gamma$, which holds for $p<1+log_{2}(\gamma/\upsilon)$. Thus, \eqref{Linearizing V 1 Agg} can be tightened using the following valid constraint:

\begin{equation}\label{Linearizing V 1 tightened p}
\Bar{V}_{p} \leq 2^{p-1}\upsilon Z_{p}, \quad   p  : p<1+log_{2}(\gamma/\upsilon)
\end{equation}

\section{Computational Results}\label{Computational Results}

In this section, we present computational results for models with the proposed constraints. Computational experiments are conducted on a Windows 10 machine with Intel Core i5 at 2.70 GHz and 8 GB of RAM. Models are coded in GAMS 33.2 and solved using CPLEX 12.10 with default settings. Instances are modified from the 90 randomly generated instances in D’Ambrosio et al. \cite{10.1007/978-3-642-20807-2_10}, which are included in QPLIB, a library of quadratic programming instances \cite{Furini2019}. The 90 instances contain 15 streams, 5–10 pools, 10 products, and 1–4 properties.

\subsection{Model $\textbf{M}^{\text{SB}}$}\label{MSB}

We first test the proposed constraints on model $\text{M}^{\text{SB}}$. Specifically, we consider the following variants:\newline

\begin{tabular}{ll}
$\text{M}^{\text{SB}}_{\text{F}}$: & Model $\text{M}^{\text{SB}}$ with \eqref{Rounded Cut 1} - \eqref{Rounded Cut 3}. \\

$\text{M}^{\text{SB}}_{\text{T}}$: & Model $\text{M}^{\text{SB}}$ with \eqref{Linearizing V 1 tightened p}. \\

$\text{M}^{\text{SB}}_{\text{FT}}$: & Model $\text{M}^{\text{SB}}$ with \eqref{Rounded Cut 1} - \eqref{Rounded Cut 3}, and \eqref{Linearizing V 1 tightened p}. \\
\end{tabular} \newline

where the proposed constraints are added for each $(j,k)$ pair.

We test our method with two different levels of discretization $n = 4$ and $n = 5$, where the finest interval for $R_{jk}$ is 1/8 and 1/16, respectively. We show the performance profiles for the instances that satisfy the following two criteria: (1) instances with at least one of the MILP models solved to optimality within 300 seconds; and (2) instances that are not solved by the slowest MILP model for that instance within five seconds.

\begin{figure}[!htb]
\includegraphics[width=\textwidth]{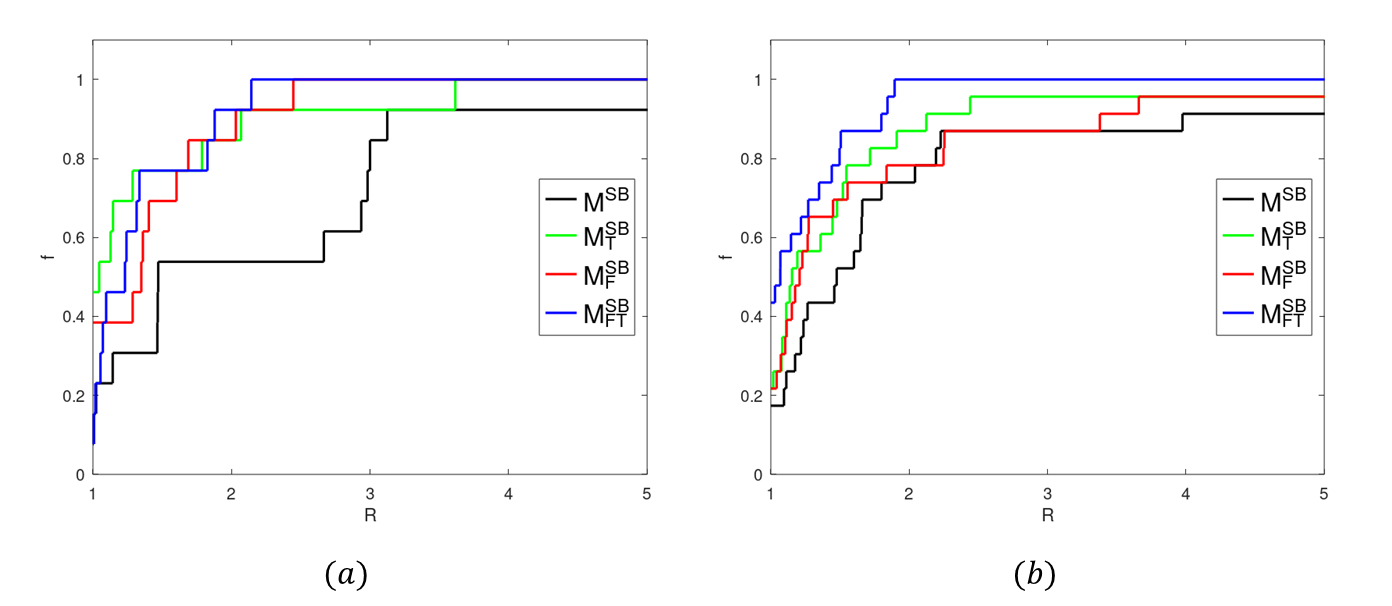}
\caption{Performance profiles for $\text{M}^{\text{SB}}$ and its variants  (a) $n = 4$.  (b) $n = 5$.} \label{PerformanceProfile_SB}
\end{figure}

The performance profile shown in Fig. \ref{PerformanceProfile_SB}(a) contains 13 instances that satisfy such criteria with $n = 4$, where the horizontal axis is the factor for performance ratio (which is defined as the solution time for a model to solve an instance over the shortest solution time among all models for the same instance), and the vertical axis is the fraction of instances \cite{Dolan2002}. We compare the optimal objective function values for these 13 instances with (1) MILP model with $n = 5$ and (2) the original NLP model in Appendix C, which demonstrates the effectiveness of the discretization method in finding a good solution to the pooling problem and further shows the value in improving the discretization-MILP models. The performance profile shown in Fig. \ref{PerformanceProfile_SB}(b) contains 23 instances that satisfy the same two criteria mentioned above with $n = 5$.

We note that the proposed constraints bring noticeable computational improvements. The model with both sets of the proposed constraints has the best performance overall. The impact of the proposed constraints is less profound with $n = 5$, when compared with $\text{M}^{\text{SB}}$, since the finer discretization intervals make the original model tighter. 


\subsection{$pq-$formulation}

We test the proposed constraints on models based on another pooling formulation, known as the $pq$-formulation \cite{Tawarmalani2002}, in which we have a nonnegative continuous variable $q_{ij} \in [0,1]$ for the proportion of stream $i$ within the total outlet flow from pool $j$, and the following nonlinear constraint:

\begin{equation}\label{nonlinear pq}
\hat{F}_{ijk} = q_{ij}\bar{F}_{jk}, \quad \forall  i, j  , k 
\end{equation}

Summing over index $k$ for the above equation we have: 

\begin{equation}\label{nonlinear pq agg}
\sum_k \hat{F}_{ijk} = q_{ij}\sum_k\bar{F}_{jk}, \quad  \forall i, j   
\end{equation}

Note that the LHS of \eqref{nonlinear pq agg} is upper bounded by the pipeline capacity between stream $i$ and pool $j$, which can be a nontrivial upper bound since the right hand side (RHS) of \eqref{nonlinear pq agg} is bounded by the capacity of pool $j$. Similar to discretizing $R_{jk}$ in $\text{M}^{\text{SB}}$, here we discretize $q_{ij}$ to obtain an MILP model, which we refer to as $\text{M}^{\text{PQ}}$. We consider variants of $\text{M}^{\text{PQ}}$ that contain the proposed constraints similar to those for $\text{M}^{\text{SB}}$, and show their performance profiles in Fig.\ref{PerformanceProfile_SB}, in which performance profiles are based on 11 (Fig.\ref{PerformanceProfile_SB}(a)) and 38 (Fig.\ref{PerformanceProfile_SB}(b)) instances that satisfy the aforementioned criteria. We observe similar improvement when adding the proposed constraints to the models based on $pq$-formulation. 

\begin{figure}[!htb]
\includegraphics[width=\textwidth]{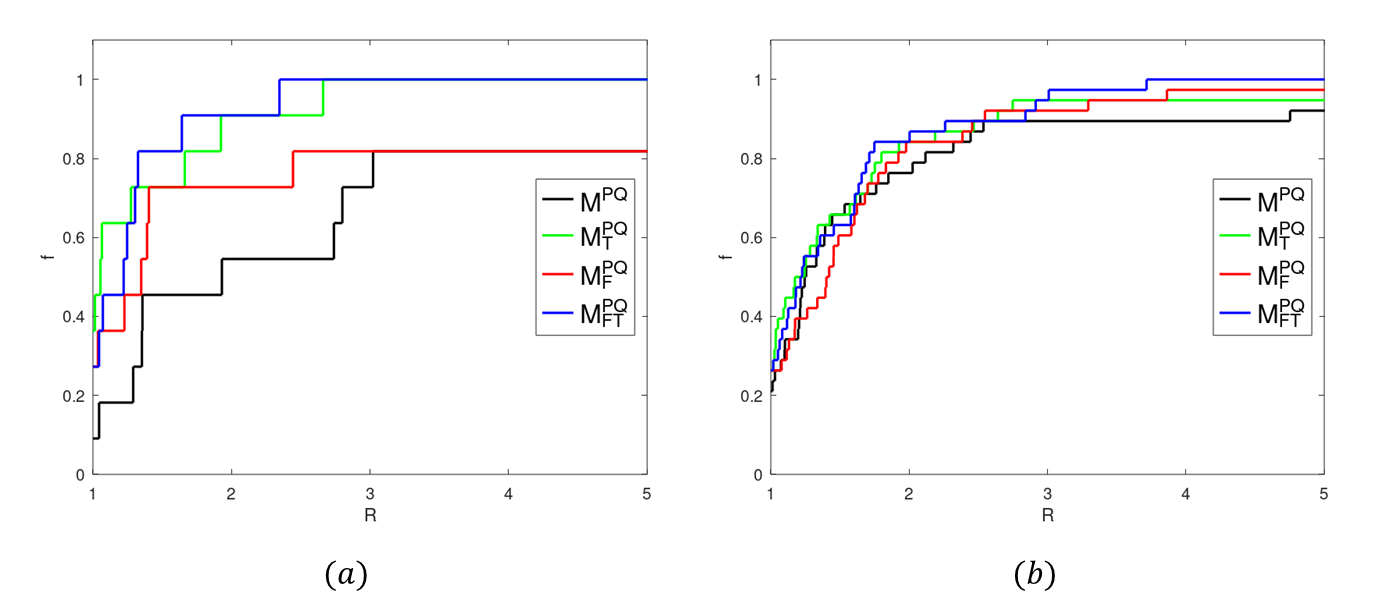}
\caption{Performance profiles for $\text{M}^{\text{PQ}}$ and its variants  (a) $n = 4$.  (b) $n = 5$.} \label{PerformanceProfile_PQ}
\end{figure}


\subsection{Alternative discretization}
We present another formulation for the discretization and demonstrate how the proposed method can be applied to it. Similar to Section 3 here we focus on a $(j,k)$ pair and drop indices $j$ and $k$. We start from the variables defined in $\text{M}^{\text{SB}}$ and explicitly define possible values that $R$ can take as $\psi_{m}$, where $m$ is the index for possible values. We also define a binary variable $Z_{m}^{'}$ and the following constraints:

\begin{equation}\label{Z_prime SOS 1}
\sum_m Z_{m}^{'} = 1   
\end{equation}

\begin{equation}\label{Z_prime SOS 1 R}
R = \sum_m \psi_{m}Z_{m}^{'} 
\end{equation}

The resulting MILP model is henceforth referred to as $\text{M}^{\text{SB-N}}$. The constraints we proposed in Section \ref{Valid Constraints} are added to $\text{M}^{\text{SB-N}}$ with minor modifications. Specifically, we now have $\mu^{-} = \text{max}_m(\psi_{m}: \psi_{m} \leq \upsilon/\gamma)$ and $\mu^{+} = \text{min}_m(\psi_{m}: \psi_{m} > \upsilon/\gamma)$, and \eqref{Rounded Cut 1} - \eqref{Rounded Cut 3} remain the same. For \eqref{Linearizing V 1 tightened p} we now have:

\begin{equation}\label{Linearizing V 1 tightened p SOS}
\Bar{V}_{m} \leq \frac{\upsilon}{\psi_{m}}Z_{p}^{'}, \quad  m  : \psi_{m} > \frac{\upsilon}{\gamma} 
\end{equation}

where $\Bar{V}_{m}$ models the bilinear term $\tilde{F}Z_{m}^{'}$. 

\begin{figure}[!htb]
\centering
\includegraphics[width=0.5\textwidth]{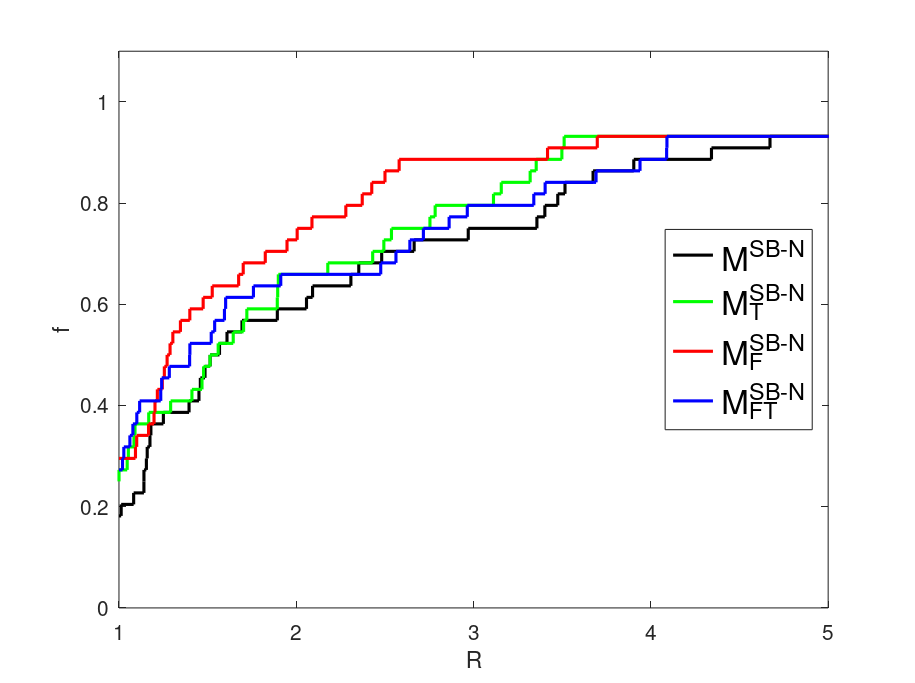}
\caption{Performance profiles for $\text{M}^{\text{SB-N}}$ and its variants.} \label{PerformanceProfile_SOS_P4}
\end{figure}

We define $\psi_{m}$ in such a way that $R$ can take the same set of values as in $\text{M}^{\text{SB}}$ with $n = 4$. We show performance profiles for variants of $\text{M}^{\text{SB-N}}$ in Fig.\ref{PerformanceProfile_SOS_P4}, which is based on 47 instances that satisfy the aforementioned criteria. We observe improvement from introducing the proposed constraints as well. 

\subsection{Comparison with lifted tangent inequalities}\label{LTI}

Valid constraints that exploit nontrivial bounds on bilinear terms with two continuous variables have been proposed \cite{Belotti2010,Belotti2011,Anstreicher2021,Chen2021}. One class of these constraints is called lifted tangent inequalities (LTI), which is involved in the convex hull of bilinear terms with nontrivial upper bounds. Here, for $\text{M}^{\text{SB}}$ and for each $(j,k)$ pair we consider the LTI in the following form given in \cite{Anstreicher2021}:

\begin{equation}\label{LTI-2}
    \frac{\upsilon}{\rho}R + \rho\tilde{F} - 2\bar{F} \leq 0
\end{equation}

where $\rho$ is a parameter. Inequality \eqref{LTI-2} is tangent to the curve $\tilde{F}R = \upsilon$ at point $(\rho, \frac{\upsilon}{\rho})$ when $\bar{F} = \upsilon$. To add \eqref{LTI-2} into the discretized model, we set the values of $\rho$ is such a way that the resulting points of tangency are feasible in $\text{M}^{\text{SB}}$. Specifically, we have:

\begin{equation}\label{LTI_p}
    \frac{\upsilon}{2^{1 - p}}R + 2^{1 - p}\tilde{F} - 2\bar{F} \leq 0, \quad \forall p \in [\![n]\!]
\end{equation}

As its name suggested, LTI can be derived from lifting the tangent line to $\tilde{F}R = \upsilon$ to include the origin. It can be strengthened for $\text{M}^{\text{SB}}$. Since only certain points on the curve are feasible in $\text{M}^{\text{SB}}$, we derive a constraint by lifting the secant line passing through two adjacent feasible points on the curve. We have:

\begin{equation}\label{LTI_p_2}
    \frac{\upsilon}{2^{2 - p} + 2^{1 - n}}R + \frac{2^{1 - p} (2^{1 - p} + 2^{1 - n})}{2^{2 - p} + 2^{1 - n}}\tilde{F} - \bar{F} \leq 0, \quad \forall p \in [\![n-1]\!]
\end{equation}

\begin{figure}[!htb]
\centering
\includegraphics[width=0.5\textwidth]{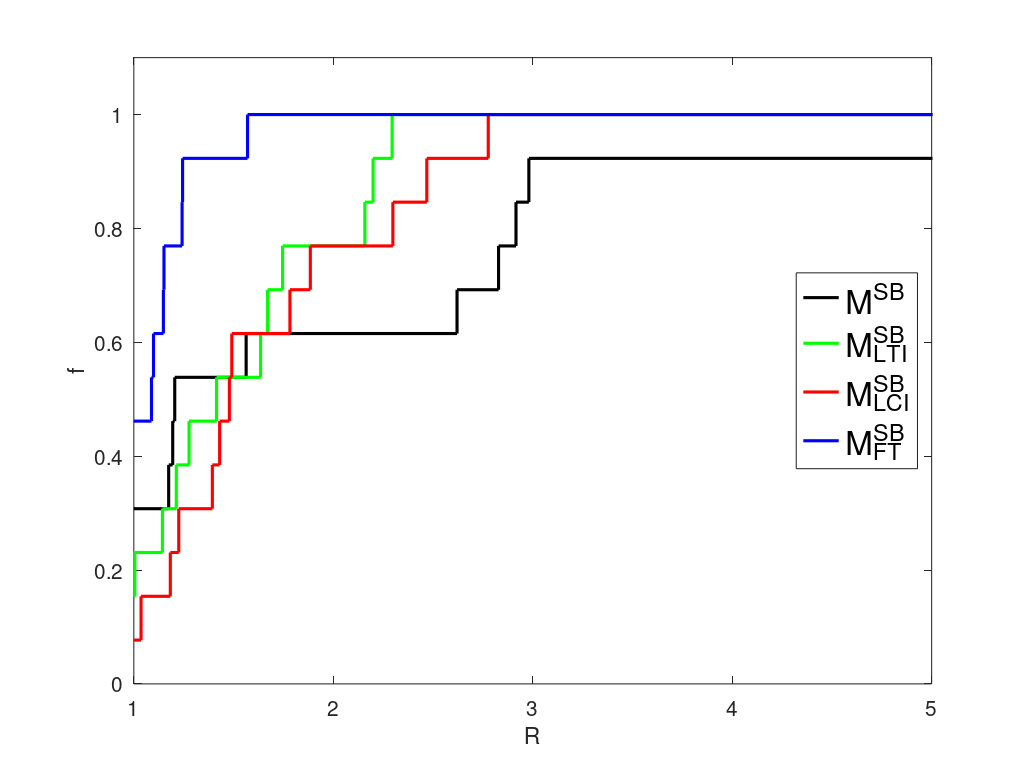}
\caption{Performance profiles for variants of $\text{M}^{\text{SB}}$ with LTI, strengthened LTI, and constraints proposed in this work.} \label{PerformanceProfile_SB_LTI}
\end{figure}

We test \eqref{LTI_p} and \eqref{LTI_p_2} on $\text{M}^{\text{SB}}$ and compare them with the constraints proposed in this work. Fig.\ref{PerformanceProfile_SB_LTI} shows the performance profiles for models with different constraints on the same set of instances tested in Section \ref{MSB} with $n = 4$, where the blue profile corresponds to the model with constraints proposed in this work, and green and red profiles correspond to models with \eqref{LTI_p} and \eqref{LTI_p_2}, respectively. We note that for the discretized model the constraints proposed in this work bring more improvement.

\section{Conclusion}\label{Conclusion}

In this paper we studied tightening methods for discretization-based MILP models for the pooling  problem, and derived valid constraints using upper bounds on bilinear terms. The proposed constraints are tested on different formulations as well as different discretization methods. Computational results demonstrate the effectiveness of our methods in terms of reducing solution time.

\section*{Acknowledgements}
The authors acknowledge financial support from the National Science Foundation under grant CBET-2026980. The authors also thank Dr. Nathan Adelgren for his comments.

\section*{Data Availability Statement}
The datasets generated during and/or analysed during the current study are available from the corresponding author on reasonable request.

\newpage

\section*{Appendices}

\subsection*{Appendix A. Comparison with the mixed-integer rounding cut}

We illustrate the difference between the constraints we proposed in Section \ref{Rounding cuts} and the well-known mixed-integer rounding cuts with the illustrative example with $\upsilon = 2.2$ shown below.

\begin{figure}[h]
\centering
\includegraphics[width=0.5\textwidth]{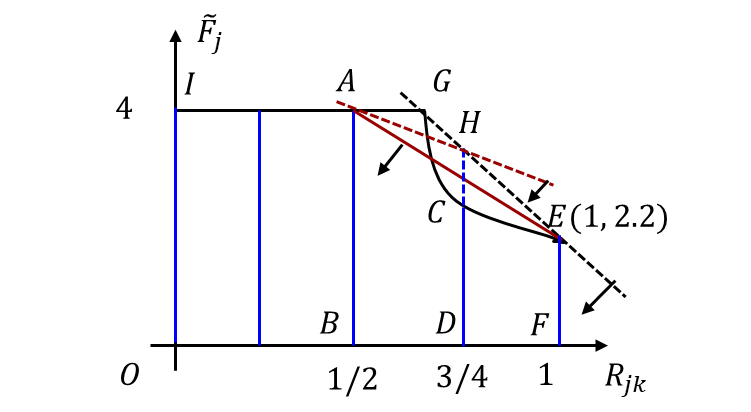}
\caption{Illustrative graph for the proposed constraint and mixed-integer rounding cut} \label{Illustrative_Appendix}
\end{figure}

To derive a mixed integer linear cut, one must start with a mixed-integer linear set. From $\textbf{S}^{*}$, the only immediately available mixed-integer linear set is obtained from using McCormick envelope to relax the constraint with bilinear term. In this example, such an approach lead to a set whose continuous relaxation has the feasible space defined by polyhedron $OIGEF$ shown in Fig.\ref{Illustrative_Appendix}, from which one can derive a mixed-integer rounding cut represented by the red dotted line passing through points $A$ and $H$. Though such a cut defines a facet for the mixed-integer linear set, we note that it is weaker than the cut that we proposed (red solid line passing through points $A$ and $E$).   

\subsection*{Appendix B. Proof of Proposition 1}

\textbf{Proposition 1.} conv($\textbf{S}^{*}$) $=\{(R,\tilde{F}) \in \mathbb{R}_{+}^{2}:  \tilde{F}\leq \gamma, R \leq 1, \eqref{Rounded Cut 1} - \eqref{Rounded Cut 3} \}$.

\begin{proof}
In this proof we consider the two cases,  $\delta \geq \epsilon$ and $\delta < \epsilon$, separately. For each case, we consider the following definition of convex hull: The convex hull of a set is the set of all convex combinations of points in it.

\textit{Case 1:} $\delta \geq \epsilon$

We define the polyhedron

\begin{equation*}
Q_1 = \{(R,\tilde{F}) \in \mathbb{R}_{+}^{2}:  \tilde{F} \leq \gamma, R \leq 1, \eqref{Rounded Cut 1} \}
\end{equation*}

We next show that (1) $Q_1$ contains all convex combinations of points in $\textbf{S}^{*}$, and (2) all points in $Q_1$ are convex combinations of points in $\textbf{S}^{*}$.

To show (1), we first note that since all points in $\textbf{S}^{*}$ are in (finitely many) line segments, all convex combinations of those points can be represented by convex combinations of the end points for those line segments. Since all end points are in the form of ($R,0$) and ($R,\text{min}(\gamma, \upsilon/R)$), they can be partitioned into three groups: (a) ($R,0$), (b) ($R,\gamma)$ with $R \leq \mu^{-}$, and (c) ($R,\upsilon/R)$ with $R \geq \mu^{+}$. It is clear that points in group (a) and (b) are in $Q_1$. We show that points in group (c) are also in $Q_1$.

For points in group (c), the inequality $R \leq 1$ is valid; To show these points also satisfy $\tilde{F} \leq \gamma$, we note that when $R = \mu^{+}$, the $\tilde{F}$-coordinate values for these points satisfy

\begin{align*}
    \frac{\upsilon}{R} = \frac{\upsilon}{2^{1-n}(\floor{\upsilon/(2^{1-n}\gamma)}+1)} \leq \frac{\upsilon}{2^{1-n}\upsilon/(2^{1-n}\gamma)} = \gamma
\end{align*}

and $\upsilon/R$ is monotonic decreasing w.r.t. $R$, thus points in group (c) satisfy $\tilde{F} \leq \gamma$.

It remains to be shown that points in group (c) also satisfy \eqref{Rounded Cut 1}. We note that these points are on the curve defined by the convex function $\tilde{F} = \upsilon/R$, and linear constraint \eqref{Rounded Cut 1} is upper bounding $\tilde{F}$. To show points in group (c) satisfy \eqref{Rounded Cut 1}, it suffices to show that the RHS of \eqref{Rounded Cut 1} is greater than or equal to $\upsilon/R$ when $R = \mu^{+}$ and when $R = 1$ (since for points in group (c) we have $R \in [\mu^{+}, 1]$).

When $R = \mu^{+}$, for the RHS of \eqref{Rounded Cut 1} we have:
\begin{align} 
\delta(\mu^{+} - 1) + \upsilon &= \frac{(\gamma - \upsilon)(\mu^{+} - 1)}{\mu^{-} - 1} + \upsilon \nonumber \\
&= \frac{(\gamma - \upsilon)(\mu^{-} + 2^{1 - n} - 1)}{\mu^{-} - 1} + \upsilon \nonumber \\
&= (\gamma - \upsilon)(1 + \frac{2^{1 - n}}{\mu^{-} - 1}) + \upsilon \nonumber \\
&= \gamma + \frac{2^{1 - n}(\gamma - \upsilon)}{\mu^{-} - 1} \label{Rounded Cut 1 RHS}
\end{align}

Since $\delta \geq \epsilon$, we have:

\begin{align}
    \frac{\gamma - \upsilon}{\mu^{-} - 1} &\geq \frac{\gamma - \upsilon/\mu^{+}}{\mu^{-} - \mu^{+}} \nonumber \\
    &= -\frac{\gamma - \upsilon/\mu^{+}}{2^{1 - n}} \label{condition delta geq epsilon}
\end{align}

Combining \eqref{Rounded Cut 1 RHS} and \eqref{condition delta geq epsilon} we have

\begin{align*}
    \delta(\mu^{+} - 1) + \upsilon &= \gamma + \frac{2^{1 - n}(\gamma - \upsilon)}{\mu^{-} - 1} \\
    & \geq \gamma - (\gamma - \upsilon/\mu^{+}) = \upsilon/\mu^{+}
\end{align*}

Thus, when $R = \mu^{+}$, the RHS of \eqref{Rounded Cut 1} is greater than or equal to $\upsilon/R$.

When $R = 1$, the RHS of \eqref{Rounded Cut 1} becomes $\upsilon = \upsilon/R$. Thus, we conclude that the RHS of \eqref{Rounded Cut 1} is greater or equal to $\upsilon/R$ when $R \in [\mu^{+},1]$, and thus points in group (c) are in $Q_1$.

We thus conclude that all end points for line segments in $\textbf{S}^{*}$ are in $Q_1$.

To show (2): all points in $Q_1$ are convex combinations of points in $\textbf{S}^{*}$, one can verify that the five extreme points of $Q_1$ are: $P_1 = (0,0),P_2 = (0,\gamma),P_3 = (\mu^{-},\gamma), P_4 = (1,\upsilon), P_5 = (1,0)$, and all of them are in $\textbf{S}^{*}$. Since all points in $Q_1$ are convex combinations of these extreme points, we conclude that all points in $Q_1$ are convex combinations of points in $\textbf{S}^{*}$.

We thus conclude that when $\delta \geq \epsilon$, conv($\textbf{S}^{*}$) = $Q_1$.

\textit{Case 2: }$\delta < \epsilon$

For this case we define the polyhedron

\begin{equation*}
Q_2 = \{(R,\tilde{F}) \in \mathbb{R}_{+}^{2}:  \tilde{F} \leq \gamma, R \leq 1, \eqref{Rounded Cut 2} - \eqref{Rounded Cut 3} \}
\end{equation*}

To show that $Q_2$ contains all convex combinations of points in $\textbf{S}^{*}$, we again consider the above three groups of points in $\textbf{S}^{*}$ defined in the previous case.

It can be shown that points in groups (a) and (b) are in $\textbf{S}^{*}$ using the steps for the previous case, we thus focus on group (c).  

For points in group (c), by inspection one can see $\tilde{F} \leq \gamma$ and $R \leq 1$ are valid. To show \eqref{Rounded Cut 2} - \eqref{Rounded Cut 3} are valid, similar to the previous case, we show that the RHS of \eqref{Rounded Cut 2} and \eqref{Rounded Cut 3} are greater than or equal to $\upsilon/R$ when $R = \mu^{+}$ and $R = 1$. 

When $R = \mu^{+}$, from both \eqref{Rounded Cut 2} and \eqref{Rounded Cut 3} we have $\tilde{F} \leq \upsilon/\mu^{+}$, thus their RHS equal to $\upsilon/R$. 
When $R = 1$, we first note that since $\delta < \epsilon$, we have:
\begin{equation*}
    \frac{\gamma - \upsilon}{\mu^{-} - 1} < \frac{\gamma - \upsilon/\mu^{+}}{\mu^{-} - \mu^{+}}
\end{equation*}

since $1 - \mu^{-} > 0$, for the RHS of \eqref{Rounded Cut 2} we have:

\begin{align*}
    \frac{\gamma - \upsilon/\mu^{+}}{\mu^{-} - \mu^{+}}(1 - \mu^{-}) + \gamma &> \frac{\gamma - \upsilon}{\mu^{-} - 1}(1 - \mu^{-}) + \gamma \\
    &=-(\gamma - \upsilon)+\gamma = \upsilon
\end{align*}

For the RHS of \eqref{Rounded Cut 3}, one can verify that when $R = 1$ it becomes $\upsilon$. 

Thus, we conclude that the RHS of \eqref{Rounded Cut 2} and \eqref{Rounded Cut 3} are greater than or equal to $\upsilon/R$ when $R \in [\mu^{+},1]$, and points in group (c) are in $Q_2$. Since points in group (a) and (b) are also in $Q_2$, we conclude that $Q_2$ contains all convex combinations of points in $\textbf{S}^{*}$.

To show that all points in $Q_2$ are convex combinations of points in $\textbf{S}_{jk}^{*}$, one can verify that the six extreme points of $Q_2$ are: $P_1 = (0,0),P_2 = (0,\gamma),P_3 = (\mu^{-},\gamma), P_4 = (1,\upsilon), P_5 = (1,0), P_6 = (\mu^{+},\upsilon/\mu^{+})$. We note that different from the previous case, here the slope of line segment $P_3P_4$, $\delta$, is smaller than the slope of line segment $P_3P_6$, $\epsilon$, thus we have a new extreme point $P_6$.  Since all points in $Q_1$ are convex combinations of these extreme points, which are also in $\textbf{S}^{*}$, we conclude that all points in $Q_1$ are convex combinations of points in $\textbf{S}^{*}$.

Combining the above results, we conclude that conv($\textbf{S}^{*}$) $= \{(R,\tilde{F}) \in \mathbb{R}_{+}^{2}:  \tilde{F} \leq \gamma, R \leq 1, \eqref{Rounded Cut 1} - \eqref{Rounded Cut 3} \}$\qed

\end{proof}

\subsection*{Appendix C. Comparison on the objective function value}

We compare the objective function obtained from the 13 instances included in Fig. \ref{PerformanceProfile_SB}(a) with (1) $\text{M}^{\text{SB}}$ with $n = 4$, (2) $\text{M}^{\text{SB}}$ with $n = 5$, and (3) the original NLP model, in Table \ref{tab:Objective function value comparison}. The NLP column shows the best known objective function value from global solver SCIP after 300 seconds (and the best possible value is shown in parentheses, if applicable). Discretization-based model lead to good solution, and finer discretization intervals can lead to better objective function value. The gap between the objective function value from MILP models and the (best possible) objective function value from NLP model, is shown in the parentheses under MILP models' objective function value. This gap is defined as $1 - \frac{\text{MILP-Obj}}{\text{BestPossible-NLP-Obj}}$.   

\newpage

\begin{table}[H]
\centering
\caption {Objective function value comparison} \label{tab:Objective function value comparison}
\begin{tabular}{lrrrl}
\cline{1-4}
\multirow{2}{*}{Instances} & \multicolumn{3}{l}{Objective Function Value and Gap}                                        &  \\ \cline{2-4}
                           & $n = 4$    & $n = 5$    & NLP                                                           &  \\ \cline{1-4}
1                          & \begin{tabular}[c]{@{}r@{}}21392.86 \\(1.79\%) \end{tabular} &  \begin{tabular}[c]{@{}r@{}} 21650\\ (0.61\%)\end{tabular}    & 21783                                                         &  \\
2                          & \begin{tabular}[c]{@{}r@{}} 21477.14 \\ (0.50\%) \end{tabular} & \begin{tabular}[c]{@{}r@{}} 21515.85 \\ (0.32\%) \end{tabular} & \begin{tabular}[c]{@{}r@{}}21560\\ (21585.63)\end{tabular}    &  \\
3                          & \begin{tabular}[c]{@{}r@{}} 33771.08 \\ (3.24\%) \end{tabular} & \begin{tabular}[c]{@{}r@{}} 33970.13\\ (2.67\%) \end{tabular} & \begin{tabular}[c]{@{}r@{}}33155.63\\ (34901.29)\end{tabular} &  \\
4                          & \begin{tabular}[c]{@{}r@{}} 30672.23 \\ (3.44\%) \end{tabular} & \begin{tabular}[c]{@{}r@{}} 30823.21 \\ (2.96\%) \end{tabular} & \begin{tabular}[c]{@{}r@{}}30371.49\\ (31764.87)\end{tabular} &  \\
5                          & \begin{tabular}[c]{@{}r@{}} 36665 \\ (0.03\%) \end{tabular}   & \begin{tabular}[c]{@{}r@{}} 36665 \\ (0.03\%) \end{tabular}   & \begin{tabular}[c]{@{}r@{}}36623.8\\ (36675)\end{tabular}     &  \\
6                          & \begin{tabular}[c]{@{}r@{}} 47030 \\ (0.51\%) \end{tabular}   & \begin{tabular}[c]{@{}r@{}} 47205 \\ (0.14\%) \end{tabular}    & \begin{tabular}[c]{@{}r@{}}47249.19\\ (47273.41)\end{tabular} &  \\
7                          & \begin{tabular}[c]{@{}r@{}} 37250\\ (0.11\%) \end{tabular}  & \begin{tabular}[c]{@{}r@{}} 37271.22 \\ (0.05\%) \end{tabular} & 37290                                                         &  \\
8                          &  \begin{tabular}[c]{@{}r@{}} 41700 \\ (-) \end{tabular}   & \begin{tabular}[c]{@{}r@{}} 41700 \\ (-) \end{tabular}   & 41700                                                         &  \\
9                          & \begin{tabular}[c]{@{}r@{}} 41220 \\ (0.24\%) \end{tabular}  & \begin{tabular}[c]{@{}r@{}} 41235 \\ (0.21\%) \end{tabular}   & \begin{tabular}[c]{@{}r@{}}41162.65\\ (41320)\end{tabular}    &  \\
10                         & \begin{tabular}[c]{@{}r@{}} 35904 \\ (0.20\%) \end{tabular}   & \begin{tabular}[c]{@{}r@{}} 35942 \\ (0.09\%) \end{tabular}   & \begin{tabular}[c]{@{}r@{}}35962.5\\ (35975.07)\end{tabular}  &  \\
11                         & \begin{tabular}[c]{@{}r@{}} 28715.18 \\ (0.46\%) \end{tabular}  & \begin{tabular}[c]{@{}r@{}} 28744 \\ (0.36\%) \end{tabular}   & \begin{tabular}[c]{@{}r@{}}28731.74\\ (28849.11)\end{tabular} &  \\
12                         & \begin{tabular}[c]{@{}r@{}} 21392.86 \\ (1.79\%) \end{tabular} & \begin{tabular}[c]{@{}r@{}} 21650 \\ (0.61\%) \end{tabular}  & 21783                                                         &  \\
13                         & \begin{tabular}[c]{@{}r@{}} 25863.7 \\ (0.67\%) \end{tabular}  & \begin{tabular}[c]{@{}r@{}} 25892 \\ (0.56\%) \end{tabular}   & \begin{tabular}[c]{@{}r@{}}25689.79\\ (26038.91)\end{tabular} &  \\ \cline{1-4}
\end{tabular}
\end{table}

\newpage

%
%
%
\bibliographystyle{splncs04}
\bibliography{DiscPoolingRef}

%


\end{document}